\documentclass[11pt]{amsart}
\usepackage{amsfonts,latexsym, color}
\usepackage{amsmath}
\usepackage{amssymb}

 \font\tenmsb=msbm10 at 12pt \font\sevenmsb=msbm7 at 8pt \font\fivemsb=msbm5 at
6pt
\newfam\msbfam
\textfont\msbfam=\tenmsb \scriptfont\msbfam=\sevenmsb \scriptscriptfont\msbfam=\fivemsb

\def\M{{\mathcal M}}
\def\R{{\mathbb R}}

\begin{document}
\newcommand{\reset}{\setcounter{equation}{0}}

\newcommand{\beq}{\begin{equation}}
\newcommand{\noi}{\noindent}
\newcommand{\eeq}{\end{equation}}
\newcommand{\dis}{\displaystyle}
\newcommand{\mint}{-\!\!\!\!\!\!\int}

\def \theequation{\arabic{section}.\arabic{equation}}

\newtheorem{thm}{Theorem}[section]
\newtheorem{lem}[thm]{Lemma}
\newtheorem{cor}[thm]{Corollary}
\newtheorem{prop}[thm]{Proposition}
\theoremstyle{definition}
\newtheorem{defn}[thm]{Definition}
\newtheorem{rem}[thm]{Remark}

\def \bx{\hspace{2.5mm}\rule{2.5mm}{2.5mm}} \def \vs{\vspace*{0.2cm}} \def
\hs{\hspace*{0.6cm}}
\def \ds{\displaystyle}
\def \p{\partial}
\def \O{\Omega}
\def \o{\omega}
\def \b{\beta}
\def \m{\mu}
\def \T{{\mathbb T}}
\def \ou{{\overline u}}
\def \ov{{\overline v}}
\def \D{\Delta}
\def \d{\delta}
\def \s{\sigma}
\def \e{\varepsilon}
\def \a{\alpha}
\def \hv{\widehat{v}}
\def \ww{{\widetilde w}}
\def \wv{{\widetilde v}}
\def \hw{{\widehat w}}
\def \wK{{\widetilde K}}
\def \l{\lambda}
\def\cqfd{%
\mbox{ }%
\nolinebreak%
\hfill%
\rule{2mm} {2mm}%
\medbreak%
\par%
}
\def \pr {\noindent {\it Proof:} }
\def \rmk {\noindent {\it Remark} }
\def \esp {\hspace{4mm}}
\def \dsp {\hspace{2mm}}
\def \ssp {\hspace{1mm}}

\title{Normal conformal metrics with prescribed $Q$-Curvature in $\R^{2n}$}
\date{}
\author{Xia Huang}
\address{School of Mathematical Sciences,  Key Laboratory of MEA (Ministry of Education) \& Shanghai Key Laboratory of PMMP,  East China Normal University, Shanghai 200241, China}
\email{xhuang@cpde.ecnu.edu.cn}
\author{Dong Ye}
\address{School of Mathematical Sciences,  Key Laboratory of MEA (Ministry of Education) \& Shanghai Key Laboratory of PMMP,  East China Normal University, Shanghai 200241, China}
 \email{dye@math.ecnu.edu.cn}
\author{Feng Zhou}
\address{School of Mathematical Sciences,  Key Laboratory of MEA (Ministry of Education) \& Shanghai Key Laboratory of PMMP,  East China Normal University, Shanghai 200241, China}
\email{fzhou@math.ecnu.edu.cn}
\date{}
\begin{abstract}
We consider the $Q$-curvature equation
\begin{equation}\label{0.1}
(-\D)^n u = K(x)e^{2nu}\quad\text{in} ~\R^{2n} \ (n \geq 2)
\end{equation}
where $K$ is a given non constant continuous function. Under mild growth control on $K$, we get a necessary condition on the total curvature $\Lambda_u$ for any normal conformal metric $g_u = e^{2u}|dx|^2$ satisfying $Q_{g_u} = K$ in $\R^{2n}$, or equivalently, solutions to \eqref{0.1} with logarithmic growth at infinity. Inversely, when $K$ is nonpositive satisfying polynomial growth control, we show the existence of normal conformal metrics with quasi optimal range of total curvature and precise asymptotic behavior at infinity. If furthermore $K$ is radial symmetric, we establish the same existence result without any growth assumption on $K$.
\end{abstract}
\maketitle

\noindent
{\small Key words: Normal conformal metrics, Nonconstant $Q$-curvature, Finite total curvature}

\noindent
{\small 2020 MSC: 53C18, 35J91, 35B08}
\vskip 0.8cm

\section {introduction}
\reset
In this paper, we consider the following polyharmonic equation
\begin{equation}\label{1.1'}
(-\D)^n u = K(x)e^{2nu}\quad \text{in}~~\R^{2n},
\end{equation}
where $n \geq 1$ and $K \ne 0$. This equation naturally arises in conformal geometry.

\smallskip
Let $(\M^2, g)$ be a Riemannian surface with the Gauss curvature $k_g$, consider a metric $g_\o = e^{2\o}g$ conformal to $g$, then the corresponding Gauss curvature $k_{g_\o}$ satisfies $-\D_{g}\o + k_g = k_{g_\o} e^{2\o}$. In particular, for $(\M^2, g) = (\O, |dx|^2)$ with $\O \subset \R^2$, it becomes $-\D\o = k_{g_\o} e^{2\o}$ in $\Omega$, called the prescribed Gauss curvature equation, whose study attracted a lot of attention.

\medskip
For general Riemannian manifold $(\M^{2n}, g)$, Graham-Jenne-Mason-Sparling \cite{GJMS} introduced higher order operators $P_g$ (called the GJMS operators) and the corresponding $Q$-curvature $Q_g$, which satisfy
$$P_g \o + Q_g = Q_{g_\o} e^{2n\o} \quad \mbox{if } g_\o = e^{2\o}g.$$
Their work generalized the construction of Branson-Orsted \cite{BO} for $(\M^4, g)$, corresponding to the Paneitz operator \cite{P}. Since then, the study of $Q$-curvature plays an important role in conformal geometry.

\medskip
One natural question is the prescription of $Q$-curvature, that is, whether we can find a conformal metric that gives the desired $Q$-curvature. Consider the flat model case with $(\M^{2n}, g_0) = (\R^{2n}, |dx|^2)$, then $P_{g_0} = (-\Delta)^n$, the problem is reduced to find solutions to \eqref{1.1'} for a given $Q$-curvature function $Q_{g_u} = K$.

\medskip
Let us begin with the simplest case $n=1$ and $K \equiv 1$. In {1850s, Liouville \cite{L1,L2}} showed for the first time that all local {real} solutions to {$-\Delta u = e^{2u}$} in {a simply connected domain} $\Omega \subset \R^2$, are given by {a locally univalent meromorphic} functions $h(z)$ with
$$u(x) = \ln\frac{2|h(z)|}{1+|h'(z)|^2}.$$
{Chou-Wan \cite{CW} and Brito-Hounie-Leite \cite{BHL} then extended Liouville's formula for non-simply connected domains of $\R^2$}. It means that there is a rich family of entire smooth solutions in $\R^2$, as we can take $h(z) = az+b$ or $e^{az+b}$ and so on, so people added the finite total curvature condition $1 \in L^1(\R^2, dv_{g_u})$ or equivalently $e^{2u} \in L^1(\R^2, dx)$, where $dx$ means the Lebesgue measure. In 1991, Chen \& Li \cite{CL} proved a famous classification result: Any smooth solution to
\begin{align}
\label{02}
-\Delta u = e^{2u} \;\;\mbox{in }\R^2, \quad \int_{\R^2} e^{2u} dx < \infty
\end{align}
is given by
\begin{align}
\label{sphere}
u_{x_0, \l}(x) = \ln\frac{2\l}{1+{\l^2}|x - x_0|^2}
\end{align}
with $x_0 \in \R^2$ and $\lambda > 0$. Geometrically, it means that any metric conformal to $|dx|^2$ in $\R^2$ with the Gauss curvature $k_{g_u} \equiv 1$ and finite total curvature is provided by a stereographic projection from the standard sphere $S^2$ into $\R^2$.

\medskip
In general, we say that the solution to \eqref{1.1'} has finite total curvature if $Ke^{2nu} \in L^1(\R^{2n}, dx)$, and we denote by $\Lambda_u$ the total $Q$-curvature associated to $g_u$, that is
$$\Lambda_u := \int_{\R^{2n}} Ke^{2nu} dx.$$

\medskip
Consider now the positive constant $Q$-curvature equation in higher dimensions.
\begin{align}
\label{03}
(-\Delta)^n u = (2n-1)! e^{2nu} \;\;\mbox{in }\R^{2n}, \quad \Lambda_u = (2n-1)! \int_{\R^{2n}} e^{2nu} dx < \infty.
\end{align}
Here $(2n-1)!$ is the $Q$-curvature of the Euclidean sphere $S^{2n}$. It is well-known that any stereographic projection from $S^{2n}$ to $\R^{2n}$ yields a conformal metric with the constant $Q$-curvature equal to $(2n-1)!$. In other words, the formula in \eqref{sphere} defines a family of solutions to \eqref{03}, and their total curvature $\Lambda_{u_{x_0, \l}} \equiv \Lambda(S^{2n})$, the total $Q$-curvature of $S^{2n}$.

\smallskip
However, the set of the conformal metrics with positive constant $Q$-curvature and finite total curvature for $n \geq 2$ becomes much more rich than the $\R^2$ case.

\smallskip
Using variational method, Chang-Chen \cite{CC} proved that for any $n \geq 2$ and any $0 < \gamma < \Lambda(S^{2n})$, there exist solutions to \eqref{1.1'} with $K \equiv (2n-1)!$ and $\Lambda_u = \gamma$. The condition $\Lambda_u < \Lambda(S^{2n})$ was not only suggested by technical reason, it is also necessary when $n =2$.
Indeed, Lin \cite{Lin} showed that any solution to \eqref{03} in $\R^4$ verifies $\Lambda_u \leq \Lambda(S^4)$, and
the equality holds if and only if the solution is spherical (see also \cite{Gu} for the geometric setting). Furthermore, for any $0 < \gamma < \Lambda(S^4)$, Wei \& Ye
\cite{WY} provided infinitely many non radial solutions to \eqref{03} in $\R^4$ with $\Lambda_u = \gamma$. This phenomenon is generalized for all $n \geq 2$ by Hyder \& Martinazzi \cite{HM}.
When $n\geq 3$, Martinazzi \cite{M} found a striking new phenomenon, he showed that there exist solutions to \eqref{03} for which the total curvatures are larger than $\Lambda(S^{2n})$.
Later, Huang-Ye \cite{HY}, Hyder \cite{H} proved that the total curvature $\Lambda_u$ can take arbitrary positive value, that is, for $n \geq 3$, $\gamma > 0$, there exist conformal metrics $g_u$ in $\R^{2n}$ such that $u$ resolves \eqref{03} and $\Lambda_u = \gamma$.

\medskip
When $K$ is a negative constant, there is no entire smooth solution to $\Delta u = e^{2u}$ in $\R^2$, see Theorem 1 in \cite{O}. Again, a completely different situation appears for $n \geq 2$. Hyder \& Martinazzi \cite{HM, H1} showed
that for any $n \geq 2$, $\gamma < 0$, and any polynomial $P$ with ${\rm deg}(P)\leq (2n-2)$ and $\lim_{|x|\to\infty} P(x) =  \infty$, there exists $u$ such that
\begin{align}
\label{03-1}
(-\Delta)^n u = -(2n-1)!e^{2nu} \quad \text{in} \; \R^{2n}, \;\; \Lambda_u = \gamma
\end{align}
 and
\begin{align}
\label{03-2}
u(x) = -P(x) + \frac{2|\gamma|}{\Lambda(S^{2n})}\ln |x| \rightarrow \ell \in \R \quad \mbox{as }\; |x|\to\infty.
\end{align}

All these results indicate that the situation is drastically changed in higher dimensional cases ($n \geq 2$), even for conformal metrics with constant $Q$-curvature.

\medskip
In \cite{CQY}, Chang-Qing-Yang used an important notion of {\sl normal metrics} for conformally flat manifolds to prove the formula of the asymptotic isoperimetric ratio. The normal metric was first introduced by Huber \cite{Hu}, he proved that in $\R^2$, the conformal metric is always normal if we have finite total Gauss curvature.
\begin{defn}
\label{normal}
A conformal metric $g_u$ in $\R^{2n}$ is said normal if it has finite total $Q$-curvature, and there is $C \in \R$ such that
\begin{align}
\label{n1}
u(x) = \frac{1}{\gamma_n} \int_{\R^{2n}}\ln\frac{|y|}{|x-y|}Q_u(y)e^{2nu}dy + C, \quad \forall\; x\in \R^{2n}
\end{align}
where
\begin{align}
\label{gamma}
\gamma_n =\frac{(2n-1)!}{2}|S^{2n}| = 2^{2n-1}(n-1)!\pi^n.
\end{align}
\end{defn}

People are interested to understand the geometric meaning of the normal metrics. Chang-Qing-Yang \cite{CQY} showed that for $n=2$, if the scalar curvature $R_{g_u}$ is nonnegative at infinity, then the metric must be normal; Wang-Wang \cite{WW} generalized this result for general $n$ and $R_{g_u}^-e^{2nu} \in L^n(\R^{2n}, dx)$; Recently, Li defined in \cite{Li1} a volume entropy for conformal metrics and showed that the metric $g_u$ is normal if and only if its volume entropy $\tau_u$ is finite.

\medskip

{When $K$ is a positive constant, Wei-Xu \cite{WX} showed that all normal solutions to \eqref{03} are given by the stereographic projection for any $n \geq 2$, i.e.~by \eqref {sphere}. If $K$ is a negative constant, Martinazzi \cite{M2} proved that there is no normal solution to \eqref{03-1} for any $n \geq 2$.}

\medskip
For general non-constant curvature cases, the prescribed Gauss curvature equation (i.e. $n = 1$ in \eqref{1.1'})
was intensively studied, see for instance \cite{Mc1}-\cite{Ni}, \cite{CL1}-\cite{CYZ} and the references therein.

\medskip
In $\R^4$, several studies were realized for the existence or nonexistence of normal conformal metrics with special $Q$-curvature; or equivalently normal solution to \eqref{1.1'} with $n=2$.
\begin{itemize}
\item For $K(x) = 1 - |x|^p$, Hyder-Martinazzi \cite{HM1} showed that a normal solution exists if and only if $p \in (0, 4)$, while for $K(x) = 1+|x|^p$, it is proved in \cite{HM1} that for any $p > 0$, there exists normal conformal metric;
\item The nonexistence results of normal conformal metric are established in \cite{JSTW, WYM} with $K(x) = (1-|x|^p)|x|^\beta$ for some $p$, $\beta$ and prescribed total curvature $\Lambda$.
\end{itemize}

\medskip
As far as we know, there are few results for the existence of normal conformal metrics with general non constant $Q$-curvature and $n \geq 2$. For $K(x) = |x|^p$ with $p > -2n \geq -2$, Hyder-Mancini-Martinazzi \cite{HMM} showed that a radial normal solution exists with $\Lambda_u = \Lambda(S^{2n})(1 + \frac{p}{2n})$.

\medskip
Here we are interested to understand normal conformal metrics with general $Q$-curvature, or more generally the solutions to \eqref{1.1'} with logarithmic control at infinity, i.e.~satisfying
\begin{align}
\label{0bis}
u(x) = O(\ln|x|) \;\;\mbox{ as } |x| \to \infty,
\end{align}
that we call also normal solution by abuse of notation.


\medskip
The following results (Theorems \ref{expan} and \ref{expan0} below) show that under mild polynomial growth control on $K$, any solution verifying \eqref{0bis} will have a precise logarithmic behavior as $|x| \to \infty$.
\begin{thm}
\label{expan}
Let $u$ be a solution to \eqref{1.1'} satisfying \eqref{0bis}. Assume that $K$ satisfies the polynomial growth condition
\begin{align}
\label{growp}
|K(x)| \leq C(|x|^\lambda + 1) \;\; \mbox{in }\; \R^{2n}, \quad \mbox{with }\; C > 0, \lambda \geq 0.
\end{align}
\begin{itemize}
\item If $K$ has constant sign at infinity, then $K(x)e^{2nu} \in L^1(\R^{2n})$, i.e.~the conformal metric $g_u = e^{2u}|dx|^2$ has a finite total curvature $\Lambda_u$.
\item If $K(x)e^{2nu} \in L^1(\R^{2n})$, there holds
\begin{align}
\label{ulog}
\lim_{|x|\to \infty}\frac{u(x)}{\ln|x|} = - \frac{\Lambda_u}{\gamma_n} = -\frac{1}{\gamma_n}\int_{\R^{2n}} K(x)e^{2nu}dx.
\end{align}
\end{itemize}
\end{thm}

A direct consequence is a lower bound for the total curvature of normal metrics $g_u$.
\begin{cor}
\label{lower}
Let $K$ have constant sign at infinity and satisfy \eqref{growp}, any normal solution of \eqref{1.1'} satisfies $$\Lambda_u \geq - \alpha_1(K)\gamma_n$$
where
\begin{align}
\label{alpha1}
\alpha_1(K): = \sup\left\{\beta \in \R \; | \; K(x)(1+ |x|)^{2n\beta} \in L^1(\R^{2n}) \right\}.
\end{align}
Hence a necessary condition to have normal conformal metric in $\R^{2n}$ with such nonpositive $Q$-curvature $K \ne 0$ is $\alpha_1(K) \geq 0$.
\end{cor}
Indeed, let $\alpha = -\frac{\Lambda_u}{\gamma_n}$. Using \eqref{ulog}, $K(x)(1 + |x|)^{2n(\alpha - \e)} \in L^1(\R^{2n})$ for any $\e > 0$, so $\alpha - \e \leq \alpha_1(K)$ hence $\alpha \leq \alpha_1(K)$.

\medskip
Moreover we get a refined asymptotic behavior of $u$ when $\Lambda_u > - \alpha_1(K)\gamma_n$.
\begin{thm}
\label{expan0}
Let $K$ satisfy \eqref{growp}. Let $u$ be a solution to \eqref{1.1'} in $\R^{2n}$ satisfying \eqref{0bis} and $Ke^{2nu} \in L^1(\R^{2n})$. If $\Lambda_u > - \alpha_1(K)\gamma_n$, then
\begin{align}
\label{uexpand}
u(x)=  -\frac{\Lambda_u}{\gamma_n} \ln |x|+u(0)+\frac{1}{\gamma_n}\int_{\R^{2n}} \ln |y| K(y) e^{2nu} dy +o(1),\quad\text{as}~|x|\to\infty.
\end{align}
\end{thm}

\medskip
Now we turn to the existence of normal conformal metrics when $K$ is nonpositive with $\alpha_1(K) > 0$. Notice that the quantity $\alpha_1(K)$ was first introduced by Cheng-Ni in \cite{CN2} for the Gauss curvature function $K$ to characterize the existence of solutions to problem \eqref{1.1'} with $n=1$.
\begin{thm}
\label{negative}
Let $n\geq 2$, $K \in C(\R^{2n})$ be nonpositive and satisfy the polynomial growth control \eqref{growp}. Assume that $\alpha_1(K)>0$, then for any $\alpha \in (0, \alpha_1)$, there exists a solution $u_{\alpha}$ to \eqref{1.1'} such that
\begin{align}\label{0}  \lim_{|x|\to \infty}\big[u(x) - \alpha\ln |x|\big] = \ell \in \R.
\end{align}
\end{thm}

\begin{rem}
Seeing Corollary \ref{lower} and $\alpha = -\frac{\Lambda_u}{\gamma_n}$, the above existence result is optimal, except for the borderline case $\alpha = \alpha_1(K)$, or equivalently $\Lambda_u = -\alpha_1(K)\gamma_n$.
\end{rem}

If furthermore $K$ is radially symmetric, we can obtain the similar existence result without any growth assumption on $K$.
\begin{thm}
\label{K<0radial}
Let $n \geq 2$ and $K$ be a nontrivial and nonpositive function in $C_{\rm rad}(\R^{2n})$. Suppose that $\alpha_1(K)>0$, then for any $\alpha \in (0, \alpha_1)$, there exists a radial solution $u_{\alpha}$ to \eqref{1.1'} satisfying \eqref{0}.
\end{thm}

The paper is organized as follows:
In Section 2, we begin with the study of the polyharmonic linear equation $\D^n u = f$ in $\R^{2n}$, we prove that with mild conditions on $f$, the solutions satisfying \eqref{ulog} have more precise asymptotic behavior at infinity. Then we use the study on the linear equation to establish Theorem \ref{expan} and \ref{expan0} in section 3. Under suitable condition on $K$, we show also a relationship between the total curvature $\Lambda_u$ and the quantity
$$\int_{\R^{2n}} (x\cdot \nabla K) e^{2nu} dx$$
via Pohozaev identity. Section 4 is addressed to the proof of our main existence result of normal conformal metrics with nonpositive $Q$-curvature. The consideration of the radial case will be given in Section 5.

\medskip
The following notations are used in the whole paper:
\begin{itemize}
\item $c_{0}(\R^m) = \big\{f \in C(\R^m), \lim_{|x|\to\infty} f(x) =0 \big\}$;
\item $c_{0, {\rm rad}}(\R^m) = c_{0}(\R^m) \cap C_{\rm rad}(\R^m)$
\end{itemize}
and
$$\dot{L}^1(\R^m) = \Big\{f \in L^1(\R^m), \int_{\R^m} f(x) dx =0\Big\}.$$
In this paper, $B_r(x)$ means the Euclidean ball centered at $x$ with radius $r > 0$, and we write $B_r$ when the center is at the origin. $C$, $C_i$ denote generic positive constants, whose dependance could be specified when necessary, and whose values can change from line to line.

\medskip
\section{Preliminary: Linear Theory}
\reset
We begin with the study of linear polyharmonic equation $\D^n u = f$ in the conformal dimension.

\begin{lem}\label{l2.2}
Let $v \in C^{2n}(\R^{2n})$ and $\D^n v = f$ in $\R^{2n}$ with $f$ having constant sign near $\infty$. If $v(x)= O({\ln |x|})$ at $\infty$, then $f \in L^1(\R^{2n})$.
\end{lem}

\noindent
{\sl Proof}. Without loss of generality, we assume that $f$ is nonnegative near $\infty$. Assume that $f \not\in L^1(\R^{2n})$, so is $\overline f$. Here and in the following, $\overline{g}$ stands for the spherical average of a generic function $g \in L^1_{loc}$, i.e.
$$
\overline{g}(r)= -\!\!\!\!\!\!\int_{\p B_r} g(x)d\sigma.
$$
Denote $v_k:=\Delta^{n-k}\overline{v}$ for $0\leq k \leq n$. Then $\Delta v_1 = \overline f$ in $\R^{2n}$, so that
\begin{align}
\label{l2}
r^{2n-1}v_1'(r)= \int^r_0 s^{2n-1} \overline f(s) ds \to \infty \quad \mbox{as } r \to \infty.
\end{align}
If $n = 1$, clearly $v_1 = \overline v$ and $\lim_{r \to \infty}|\overline v(r)|/\ln r =\infty$, which is impossible.

\medskip
From now on, let $n\geq 2$. By \eqref{l2}, $v_1$ is nondecreasing near $\infty$ and $\lim_{r \to \infty} {v_1}(r)=: \ell_1 \in \R\cup{\{\infty}\}$ exists. If $\ell_1\neq 0$, then $\overline{v}(r)\sim \frac{\ell_1}{C_0} r^{2n-2}$ as $r \to \infty$, which contradicts the assumption $v(x)=O(\ln |x|)$ at infinity, hence $\ell_1 = 0$.

\smallskip
Using again \eqref{l2}, given any $M > 0$, $ v_1'(r) \geq Mr^{1-2n}$ for $r$ large, therefore $$-v_1(r) = \int_0^\infty v_1'(s)ds \geq \frac{M}{2(n-1)}r^{2-2n}  \quad \mbox{for $r$ large}.$$
By induction, we can claim that for any $1 \leq k \leq n-1$, $\lim_{r\to\infty} v_k(r)=0$, and there holds $(-1)^k v_k\geq MC_k r^{2(k-n)}$ for $r$ large enough, where $C_k$ are positive constants depending only on $k$ and $n$.

\medskip
So we get $(-1)^{n-1} \D \overline v\geq CM r^{-2}$ for $r$ large enough, therefore $(-1)^{n-1}\overline v'\geq -CMr^{-1}$ for large $r$. As $M > 0$ is arbitrary, $\lim_{r \to \infty}|\overline v(r)|/\ln r =\infty$, which yields a contradiction. Hence $f$ must belong to $L^1(\R^{2n})$.
\qed

\medskip
For the existence issue to $(-\D)^n v = f$ in $\R^{2n}$ verifying \eqref{0bis}, we will use a polynomial growth condition
\begin{align}\label{grow}
|f(x)| \leq C(|x|^\lambda + 1) \quad \text{with}~\lambda \geq  0.
\end{align}

\begin{thm}
\label{log}
Let $n\geq 1$. Suppose that $f \in L^1(\R^{2n})$ satisfies \eqref{grow} in $\R^{2n}$. Then there exists a continuous function $v$ satisfying $(-\D)^n v = f$ in $\R^{2n}$ and \eqref{0bis}. Moreover, $v$ is unique up to a constant, and
\begin{align}
\label{vlog}
\lim_{|x|\to \infty}\frac{v(x)}{\ln|x|} = -\frac{1}{\gamma_n}\int_{\R^{2n}} f(x)dx.
\end{align}
\end{thm}

\noindent
{\sl Proof}. Denote
\begin{align}
v(x) = -\frac{1}{\gamma_n}\int_{\R^{2n}}\ln\frac{|x-y|}{|y|+1} f(y)dy, \quad \forall\; x \in \R^{2n}.
\end{align}
Clearly $v \in C(\R^{2n})$ is well defined and $(-\D)^n v = f$ in $\R^{2n}$. Then
\begin{align*}
\frac{\gamma_n v(x)}{\ln|x|} + \int_{\R^{2n}} f(x)dx &= \int_{\R^{2n}}\frac{- \ln|x-y| + \ln(|y|+1) +\ln|x|}{\ln|x|} f(y)dy\\
& =: \sum_{1\leq i \leq 3} \int_{T_i} P(x, y) f(y)dy\\
& =: I_1 + I_2 + I_3,
\end{align*}
where
\begin{align*}
T_1 = B_1(x),\quad T_2 = B_R, \quad T_3 = \R^{2n}\setminus (B_1(x)\cup B_R)
\end{align*}
with $R \geq 1$, $|x| \geq R + 2$.

\medskip
First, $\lim_{|x|\to\infty} I_2 =0$ because $\|P(x, y)\|_{L^\infty(T_2)} \to 0$ if $|x|\to \infty$. In $T_1$, as $|x| \geq 3$,
\begin{align*}
|I_1| \leq C\int_{B_1(x)}|f(y)|dy + \frac{1}{\ln|x|}\int_{B_1(x)} -|f(y)|\ln|x-y|dy.
\end{align*}
By \eqref{grow} and $f \in L^1(\R^{2n})$,
\begin{align*}
& \int_{B_1(x)} -|f(y)|\ln|x-y|dy\\
= & \; \int_{|x-y| \leq |x|^{-(\lambda+1)}}-|f(y)|\ln|x-y|dy + \int_{|x|^{-(\lambda+1)} < |x-y| < 1} -|f(y)|\ln|x-y|dy\\
\leq & \; C|x|^\lambda \int_{|x-y| \leq |x|^{-(\lambda+1)}}-\ln|x-y|dx + (\lambda+1)\ln|x|\int_{|x|^{-(\lambda+1)} < |x-y| < 1} |f(y)|dy\\
\leq & \; C|x|^{(1-2n)\lambda-2n}\ln|x| + o(\ln|x|).
\end{align*}
We see that $\lim_{|x|\to\infty} I_1 =0$. In $T_3$,
there holds
\begin{align*}
|P(x, y)| \leq C, \quad \mbox{since } \frac{1}{|x|} \leq \frac{|y|+1}{|x-y|} \leq \frac{|x-y| + |x|+1}{|x-y|}\leq 2 + |x|.
\end{align*}
Therefore
\begin{align*}
|I_3| \leq C\int_{\R^{2n}\setminus B_R} |f(y)|dy \longrightarrow 0 \;\;\mbox{as } R \to \infty.
\end{align*}
We conclude then
\begin{align*}
\lim_{|x|\to\infty}\left[\frac{\gamma_n v(x)}{\ln|x|} + \int_{\R^{2n}} f(x)dx\right] = \lim_{|x|\to\infty}(I_1+I_2+I_3) = 0.
\end{align*} If there is another solution $w$ satisfying \eqref{0bis}, then $\D^n(v-w) =0$ and $(v-w) = O(\ln|x|)$ at $\infty$. Standard elliptic theory yields then $v - w$ is
constant in $\R^{2n}$. \qed

\medskip
When strengthening a little bit the integrability of $f$, we get a more precise asymptotic behavior for the above solution $v$.
\begin{thm}
\label{linear}
Let $n\geq 1$. Assume that $f$ satisfies \eqref{grow} in $\R^{2n}$ and $f(x)\ln(2 + |x|) \in L^1(\R^{2n})$. Then there exists a unique continuous function $v$ such that $(-\D)^n v = f$ in $\R^{2n}$ and
\begin{align}
\label{vlog0}
\lim_{|x|\to\infty} \left[v(x)  + \frac{\ln|x|}{\gamma_n}\int_{\R^{2n}} f(y)dy\right] = 0.
\end{align}
If moreover $f \in \dot{L}^1(\R^{2n})$, there exists a unique $v \in c_0(\R^{2n})$ such that $(-\D)^n v = f$ in $\R^{2n}$.
\end{thm}

\noindent
{\sl Proof}. As $f(x)\ln(2 + |x|) \in L^1(\R^{2n})$ and $f$ is locally bounded, define now
\begin{align}
\label{vlog1}
v(x) = -\frac{1}{\gamma_n}\int_{\R^{2n}}f(y)\ln|x-y| dy, \quad \forall\; x \in \R^{2n}.
\end{align}
Then $v \in C(\R^{2n})$ is well defined and $(-\D)^n v = f$ in $\R^{2n}$. So
\begin{align*}
\gamma_n v(x)  + \ln|x|\int_{\R^{2n}} f(y)dy &= \int_{\R^{2n}}\ln\frac{|x|}{|x-y|} f(y)dy\\
& =: \sum_{1\leq i \leq 4} \int_{T_i} \ln\frac{|x|}{|x-y|}f(y)dy\\
& =: \sum_{1\leq i \leq 4} I_i
\end{align*}
where
\begin{align*}
T_1 =\left\{y \in \R^{2n}, |y - x |\leq |x|^{-(\lambda+1)}\right\}
\end{align*}
and
\begin{align*}
T_2 = B_R, \quad T_3 = B_\frac{|x|}{2}(x)\setminus T_1, \quad T_4 = \R^{2n}\setminus ( B_\frac{|x|}{2}(x)\cup T_2)
\end{align*}
with $R \geq 1$, $|x| \geq 2R + 2$.

\medskip
Obviously, $\lim_{|x|\to\infty} I_2 =0$ since $f \in L^1(\R^{2n})$. For $I_1$, we have
\begin{align*}
|I_1| &\leq \int_{T_1} |f(y)|\ln|x|dy - \int_{T_1} |f(y)|\ln|x-y|dy \\
& \leq C\int_{T_1} |f(y)|\ln|y|dy - C|x|^\lambda \int_{T_1} \ln|x-y|dy\\
& \leq C\int_{T_1} |f(y)|\ln|y|dy + C|x|^{(1-2n)\lambda -2n}\ln|x|.
\end{align*}
Moreover, there holds
\begin{align*}
2 \leq \frac{|x|}{|x-y|} \leq |x|^{\lambda+2} \leq 2^{\lambda+2}|y|^{\lambda+2} \quad \mbox{in }\; T_3
\end{align*}
and
\begin{align*}
\frac{1}{1 + |y|}\leq \frac{|x|}{|x-y|} \leq 2  \quad \mbox{in }\; T_4,
\end{align*}
since $|x - y| \leq |x| + |y| \leq |x|(1 + |y|)$ in $T_4$. Therefore,
\begin{align*}
|I_i|\leq C\int_{T_i}|f(y)|\ln(|y| + 2)dy, \quad \mbox{for } i = 3, 4.
\end{align*}
Combining all these estimates, we deduce that
\begin{align*}
\limsup_{|x|\to \infty}\Big|\gamma_n v(x)  + \ln|x|\int_{\R^{2n}} f(y)dy\Big| \leq C\int_{\R^{2n}\setminus B_R}|f(y)|\ln(|y| + 2)dy.
\end{align*}
Taking $R \to \infty$, \eqref{vlog0} holds true. When $f \in \dot{L}^1(\R^{2n})$, \eqref{vlog0} means just $v \in c_0(\R^{2n})$.\qed

\begin{rem}
As we used only the kernel to prove Theorems \ref{log} and \ref{linear}, similar results remain valid for the corresponding linear equation in odd dimensions, i.e.~for $(-\D)^{\frac{m}{2}} v=f$ in $\R^m$.
\end{rem}

\medskip
The following result is not directly used to study the $Q$-curvature equation \eqref{1.1'}, but some ideas of the approach will be useful for us.
\begin{thm}\label{expan1}
Suppose that $f$ satisfies \eqref{grow} and $f(x)(1+|x|)^{\epsilon_0} \in L^1(\R^{2n})$ for some $\epsilon_0>0$. Then there exists a unique function $v \in C(\R^{2n})$ satisfying $\D^n v = f $ in $\R^{2n}$ and
\begin{align}
\label{vlog2}
  v(x)= -\frac{\ln|x|}{\gamma_n}\int_{\R^{2n}} f(x)dx +O(|x|^{-\delta})  \quad\text{with $\delta > 0$, as $|x| \to \infty$}.
\end{align}
\end{thm}

\noindent
Proof. Let $v$ be defined by \eqref{vlog1} and
$$\ww(x) = v(x)+\frac{\ln |x|}{\gamma_n}\int_{\R^{2n}} f(x)dx. $$
Denote by $w$ the Kelvin transformation of $\ww$ with the respect to the unite disk, i.e. $w(x) = \ww(x/|x|^2)$, then
$$
\D^n w=|x|^{-4n} f\left(\frac{x}{|x|^2}\right) =: h\quad\text{in}~B_1\setminus{\{0}\},
$$
By the assumption on $f$, we have $h(x)|x|^{-\epsilon_0} \in L^1(B_1)$. Fix $\epsilon>0$ such that $(4n+\lambda)\epsilon = \epsilon_0$, using \eqref{grow}, there holds
$$
|h(x)|^{1+\epsilon} = |h(x)| |x|^{-\epsilon_0} \left(|h(x)||x|^\frac{\epsilon_0}{\epsilon}\right)^\epsilon \leq |h(x)| |x|^{-\epsilon_0} \left(C|x|^{-4n-\lambda+\epsilon_0/\epsilon}\right)^\epsilon\leq C_\epsilon |h(x)| |x|^{-\epsilon_0}.
$$
So $h\in L^{1+\epsilon}(B_1)$. Let $\xi$ be the solution of
$$
\D^n \xi=h\quad\text{in}~B_1, \quad\quad \D^i \xi|_{\partial B_1}=0,~~0\leq i\leq n-1.
$$
By classical elliptic theory and Sobolev embedding, $\xi\in C^{0,\delta}(B_1)$ for some $\delta>0$. As $w\in L^\infty(B_1)$ is bounded by Theorem \ref{linear}
$$
-\D^n\zeta=0\quad \text{and}\quad \zeta\geq 0 \quad \text{in}~B_1\setminus{\{0}\}.
$$
Applying Theorem 3.1 in \cite{GMT}, we get
$$
\zeta(x) =H(x) + \sum_{|\beta|\leq 2n-2} a_\beta D^\beta \Phi(x) \quad \text{in}~B_{\frac{1}{2}}\setminus{\{0}\},
$$
where $\beta$ are multi-index, $a_\beta \in \R$, $\Phi(x) = -\ln|x|$ and $H\in C^{\infty}(B_1)$ is a polyharmonic function. Remark that $\zeta \in L^\infty(B_1)$, we must have $a_\beta\equiv 0$ for any $\beta$. Hence,
$$w = H + \xi - C_1 \;\; \in C^{0,\gamma}(B_\frac{1}{2}).$$
Coming back to $\ww$ who tends to $0$ at infinity, we get the claimed expansion. \qed

\medskip
\section{Asymptotic behavior and Pohozaev identity}
\reset
\subsection{Proof of Theorem \ref{expan} and \ref{expan0}}
The results are direct consequences of the previous section.

\medskip
Let $u$ be a solution to \eqref{1.1'} verifying \eqref{0bis} and $K$ satisfy \eqref{growp}. Clearly $K(x)e^{2nu}$ verifies also the polynomial growth condition \eqref{growp} with a different $\lambda$. Suppose now $K$ is of constant sign at infinity, so is $K(x)e^{2nu}$. Lemma \ref{l2.2} implies that $K(x)e^{2nu} \in L^1(\R^{2n})$, then \eqref{ulog} holds true.
The proof of Theorem \ref{expan} is completed.

\medskip
If we assume that the total curvature $\Lambda_u > -\alpha_1(K)\gamma_n$, using \eqref{ulog}, there exists $\epsilon_0 > 0$ such that $K(x)e^{2nu} \in L^1(\R^{2n}, (1 + |x|)^{\epsilon_0} dx)$. Let $v$ be given by Theorem \ref{linear}, or more precisely by \eqref{vlog1} with $f = Ke^{2nu}$. As $u-v = O(\ln|x|)$ is polyharmonic in $\R^{2n}$, there holds $u - v \equiv \ell$ in $\R^{2n}$, hence
$$
u(x)=u(0)-\frac{1}{\gamma_n} \int_{\R^{2n}}(\ln |x-y|-\ln |y|)K(y)e^{2nu} dy,
$$
Applying \eqref{vlog0}, Theorem \ref{expan0} holds true.\qed

\begin{rem}
If $\Lambda_u > -\alpha_1(K)\gamma_n$ in Theorem \ref{expan0}, we could also apply Theorem \ref{expan1} to get more precise behavior for the solution $u$. Moreover, Theorem \ref{expan} is in fact valid for the $Q$-curvature equation $(-\Delta)^\frac{m}{2} u = K(x) e^{mu}$ in $\R^m$ for any dimension $m \geq 1$.
\end{rem}

\subsection{A Pohozaev Identity}
Under more precise growth condition like
\begin{align}
\label{growp1}
C_1 |x|^\lambda \leq |K(x)| \leq C_2|x|^\lambda, \quad \text{for}~|x|~\to \infty,
\end{align}
with $C_1, C_2 > 0$ and $\lambda \in\R$, we have the following equality which is a consequence of the Pohozaev identity.

\begin{thm}\label{pohozaev}
Let $K \in C^1(\R^{2n})$ satisfy \eqref{growp1}. Let $u$ be a solution of \eqref{1.1'} verifying \eqref{0bis}. Suppose that $K(x)e^{2nu} \in L^1(\R^{2n})$ and $\Lambda_u > -\alpha_1(K)\gamma_n$. Assume moreover either $x\cdot \nabla K$ is of constant sign at infinity or $\alpha_1(x\cdot \nabla K) \geq \alpha_1(K)$, then
\begin{align}
\label{P1}
\int_{\R^{2n}} (x\cdot \nabla K) e^{2nu} dx = \frac{n\Lambda_u}{\gamma_n}( \Lambda_u - 2\gamma_n).\end{align}
\end{thm}

\noindent
{\sl Proof}. By \eqref{growp1}, obviously $\alpha_1(K) = -1 - \frac{\lambda}{2n}$. As $\Lambda_u > -\alpha_1(K)\gamma_n$, by the proof of Theorem \ref{expan0}, there holds
$$
 u(x) = -\frac{1}{\gamma_n}\int_{\R^{2n}} K e^{2nu} \ln |x-y| dy + C, \quad \forall\; x \in \R^{2n}.
$$

For $k \geq 1$, denote
$$D_k w := \left\{ \begin{array}{ll} \D^\frac{k}{2}w &\mbox{if $k$ is even};\\
\ds \p_r\big(\D^\frac{k-1}{2}w\big) = \frac{x}{|x|}\cdot\nabla\big(\D^\frac{k-1}{2}w\big) &\mbox{if $k$ is odd}.
\end{array} \right.$$
We claim that
\begin{align}
\label{Dku}
\forall\; 1\leq k\leq 2n-1, \quad \lim_{|x| \to +\infty}|x|^k D_k u(x) = A_k\alpha \;\; \mbox{with } \alpha = -\frac{\Lambda_u}{\gamma_n}
\end{align}
and
$$
A_k= (-1)^{\left[\frac{k-1}{2}\right]}2^{k-1} \left[\frac{k-1}{2}\right]!\frac{(n-1)!}{\left(n-1 - \left[\frac{k}{2}\right]\right)!}.
$$
Here
$[t]$ means the integer part of $t \in \R$.

\medskip
Indeed, $D_k\ln|x| = A_k|x|^{-k}$ in $\R^{2n}\backslash\{0\}$ for any $1\leq k\leq 2n-1$. There holds then
\begin{align*}
|x|^k D_k u(x) - A_k\alpha & = \frac{1}{\gamma_n}\int_{\R^{2n}}|x|^k\Big(D_k \ln|x|-D_k\ln|x-y|\Big) K e^{2nu} dy\\
& =: \sum_{i = 1}^3 \frac{1}{\gamma_n} \int_{T_i} P_1(x,y) K e^{2nu} dy\\
& =: I_1 + I_2 + I_3
\end{align*}
where
\begin{align*}
T_1 =\left\{y \in \R^{2n}, |y - x |\leq \frac{|x|}{2}\right\},\;\; T_2 = B_R, \;\; T_3 = \R^{2n}\setminus ( T_1\cup T_2)
\end{align*}
with $R\geq 1$, $|x|\geq 2R+2$.

\medskip
Thanks to the mean value theorem, for any $y \in T_2$
$$
\Big|D_k\ln|x|- D_k\ln |x-y|\Big| \leq |y|\times \sup_{ B_R(x)} \big|\nabla^{k+1}\ln |z|\big|\leq \frac{C}{|x|^{k+1}},
$$
so we get
\begin{align*}
|I_2| \leq \frac{C}{|x|}\int _{\R^{2n}}|K|e^{2nu} dx \to 0 \quad \mbox{ as } \; |x|\to \infty.
\end{align*}
In $T_3$, there holds $|Q(x,y)|\leq C$, then
$$
|I_3|\leq C\int_{\R^{2n}\setminus{B_R}} |K| e^{2nu} dx\to 0, \quad \text{as}~R\to \infty.
$$
Finally, for $|x|$ large, we have $|K(y)|e^{2nu}\leq C|x|^{2n \alpha+\lambda}$ in $T_1$. Therefore, let $\alpha = -\frac{\Lambda_u}{\gamma_n}$,
\begin{align*}
|I_1|& \leq C\int_{T_1} \left(1 + \frac{|x|^k}{|x-y|^k}\right) |K| e^{2nu} dy\\
& \leq C\int_{T_1} \left(1 + \frac{|x|^k}{|x-y|^k}\right) |x|^{2n\alpha+\lambda} dy \\
& \leq C|x|^{2n\alpha+\lambda+2n}.
\end{align*}
So $\lim_{|x|\to\infty} I_1 =0$ because $\alpha < \alpha_1(K) = -1 - \frac{\lambda}{2n}$. The claim \eqref{Dku} holds true. Similarly, we have also
\begin{align}
\label{Dku2}
 |x|^k D_k(x\cdot \nabla u)\to 0 \quad \mbox{as $|x|\to \infty$}, \quad \forall\; 1\leq k\leq 2n-1.
\end{align}

\medskip
Applying the Pohozaev identity for $(-\Delta)^n u= K(x) e^{2nu}$ in $B_R$ for $R > 0$, there holds,
\begin{align}
\label{Poho1}
\begin{split}
&\frac{1}{2n}\int_{B_R}(x\cdot \nabla K) e^{2nu} dx + \int_{B_R} K e^{2nu} dx + \frac{R}{2}\int_{\partial{B_R}}(D_n u)^2d\sigma\\=&\; \frac{R}{2n}\int_{\partial{B_R}} K e^{2nu} d\sigma - \sum_{j=1}^{n-1}\int_{\partial{B_R}} (-1)^{n+j}(D_{2n-1-j} u)D_j(x\cdot\nabla u)d\sigma\\
& \; + (-1)^{n+1}R\int_{\partial{B_R}} (D_{2n-1} u)(D_1 u)d\sigma.
\end{split}
\end{align}
For the last term, we used $(x\cdot \nabla u) = R (D_1u)$ on $\partial B_R$. By \eqref{Dku}-\eqref{Dku2}, \eqref{growp1} and $\alpha < \alpha_1(K)$, the right hand side goes to $(-1)^{n+1} A_{2n-1} A_1 \alpha^2 |S^{2n-1}|$ as $R\to \infty$.

\medskip
If $(x\cdot \nabla K)$ is of constant sign at infinity or $(x\cdot \nabla K) e^{2nu} \in L^1(\R^{2n})$ (when $\alpha_1(x\cdot \nabla K) \geq \alpha_1(K)$), tending $R \to \infty$ in \eqref{Poho1}, using \eqref{Dku} and $|S^m| = (m+1)\pi^\frac{m+1}{2}/\Gamma\left(\frac{m+3}{2}\right)$,  we deduce that
\begin{align}\label{p}
\frac{1}{2n}\int_{\R^{2n}} (x\cdot \nabla K) e^{2nu} dx & = \gamma_n\alpha - \frac{|S^{2n-1}|}{2}A_n^2\alpha^2 + (-1)^{n+1} A_{2n-1} A_1 \alpha^2 |S^{2n-1}|\\\nonumber
& = \gamma_n\alpha + \frac{1}{2}\frac{2\pi^n}{(n-1)!}\left[2^{n-1}(n-1)!\right]^2\alpha^2\\\nonumber
& = \frac{\gamma_n}{2}\alpha(\alpha+2),
\end{align}
{with $\alpha=-\frac{\Lambda_u}{\gamma_n}$.} \qed

\begin{rem}
{For $K(x)=|x|^p$, as the left hand side of \eqref{P1} equals to $p\Lambda_u$, we see that any normal solution to \eqref{1.1'} satisfies
\begin{align*}
\alpha=-2-\frac{p}{n},\quad \mbox{or equivalently } \;\Lambda_u= \Lambda(S^{2n})\Big(1+\frac{p}{2n}\Big).
\end{align*}
This was proved in \cite{HMM} for radial normal solutions (see also \cite{GHYZ} with $n=2$).}\end{rem}


\medskip
\section{Normal conformal metrics for negative $Q$-curvature}
Here we will prove the existence of normal solutions to \eqref{1.1'} for $K \leq 0$, i.e. Theorem \ref{negative}. Fix $u_0 \in C^\infty(\R^{2n})$ such that $u_0(x) = \ln|x|$ for any $|x|\geq \frac{1}{2}$. Clearly, $(-\D)^n u_0\equiv 0$ in $\R^{2n}\setminus B_\frac{1}{2}$ and
\begin{align*}
\int_{\R^{2n}}(-\Delta)^n u_0 dx = - \gamma_n.
\end{align*}
 We would like to find a solution $u$ to \eqref{1.1'} of the form $$u=w+\alpha u_0+\frac{1}{2n}c_w$$
 with $w\in c_0(\R^{2n})$ and a suitable constant $c_w$. Indeed, a direct computation shows that to get $u$ solution of \eqref{1.1'} is equivalent to find $w \in c_0(\R^{2n})$, $c_w\in\R$ such that
\begin{align}
\label{new1.1}
(-\Delta)^n w = K(x) e^{2n\alpha u_0 + 2nw+c_w} - \alpha(-\Delta)^n u_0 \;\; \mbox{in $\R^{2n}$}.
\end{align}
To have such $w$, we will use Leray-Schauder's fixed point argument, see Theorem 11.3 in \cite{GT}.
\begin{lem}
\label{LS} Let $T$ be a compact mapping of a Banach space $X$ into itself, and suppose that there exists a constant $M_0$ such that $\|z\|_X \leq M_0$
for all $z\in X$ satisfying $z=tT(z)$ with some $t \in (0, 1]$. Then $T$ has a fixed point.
\end{lem}

\medskip\noindent
{\sl Proof of Theorem \ref{negative}}. We will proceed by several steps.

\medskip
{\sl Step 1: The general setting.} Fix $\alpha\in(0,\alpha_1)$. For any $w \in c_0(\R^{2n})$, there exists a unique $c_w \in \R$ such that
\begin{align}
\label{cw}
\int_{\R^{2n}} K(x) e^{2n\alpha u_0 + 2nw+c_w} dx = -\alpha \gamma_n.
\end{align}
 Indeed, $c_w$ is given by
\begin{align}\label{c_w}
c_w=\ln (\alpha \gamma_n)-\ln \left\{\int_{\R^{2n}}- K(x) e^{2n\alpha u_0 + 2nw} dx\right\}.
\end{align}
Now we can apply Theorem \ref{linear} to get a solution $v \in c_0(\R^{2n})$ such that
\begin{align*}
(-\Delta)^n v = K(x) e^{2n\alpha u_0 + 2nw+c_w} - \alpha(-\Delta)^n u_0 =: \Psi(w)\quad \mbox{in $\R^{2n}$}.
\end{align*}
Since $\alpha < \alpha_1(K)$, there exists $\e_0 > 0$ such that
\begin{align*}
K(x) e^{2n\alpha u_0 + 2nw+c_w}(1+|x|)^{\e_0} \in L^1(\R^{2n}),
\end{align*}
hence $\Psi(w) \in \dot{L^1}(\R^{2n})\cap L^1(\R^{2n}, \ln(2+|x|)dx)$. Denote $\T(w) = v$, so $\T$ is a mapping from $c_0(\R^{2n})$ into itself and it's not difficult to see that $\T$ is continuous.

\medskip
{\sl Step 2: Compactness of the operator $\T$.} Recall that $c_0(\R^{2n})$ is a Banach space with the norm $\|\cdot\|_\infty$. Consider the ball $\Sigma := \{\|w\|_\infty \leq M\}$ in $c_0(\R^{2n})$ with $M > 0$, we shall prove that $\T(\Sigma)$ is relatively compact in $c_0(\R^{2n})$.

\smallskip
Using \eqref{c_w}, we see from that $c_w$ are uniformly bounded for $w \in \Sigma$. Therefore the family $\{\Psi(w), w \in \Sigma\}$ is uniformly integrable in $L^1(\R^{2n}, \ln(2+|x|)dx)$ and verifies the growth condition \eqref{growp} with some common constants $C$ and $\lambda$.

\smallskip
By the proof of Theorem \ref{linear}, for any $\e > 0$, there exists $R > 0$ depending only on $M$ and $\e$, such that
\begin{align}
\label{epsilon}
|\T(w)(x)|\leq \e, \quad \forall\; |x| \geq R,\; w \in \Sigma.
\end{align}
Furthermore, we have the integral representation of $\T(w)$,
\begin{align}
\T(w)=-\frac{1}{\gamma_n}\int_{\R^{2n}} \ln |x-y| K(x) e^{2n\alpha u_0 + 2nw+c_w}-\alpha u_0(x).
\end{align}
There holds then
$$|\Delta(\T(w)+\alpha u_0)(x)|\leq C_n\int_{\R^{2n}} \frac{| K(y)|}{|x-y|^2} e^{2n\alpha u_0 + 2nw+c_w} dy,
$$
which yields readily that
\begin{align*}
|\T(w)(x)| + |\D\T(w)(x)| \leq C_{M, R}, \quad \forall\; |x| \leq R+1, \; w \in \Sigma.
\end{align*}
By elliptic theory, $\{\T(w)|_{\overline B_R},\; w \in \Sigma\}$ is relatively compact in $C(\overline B_R)$. Combining with \eqref{epsilon}, the family $\{\T(w),\; w \in \Sigma\}$ is precompact, hence relatively compact in $c_0(\R^{2n})$. So the operator $\T$ is compact.

\medskip
{\sl Step 3: Uniform local estimate.} Now we consider the family of eventual fixed points of $ t\T(w)$ with $t \in (0, 1]$. Assume that $w=t\T(w)$, $t\in(0,1]$ and $w\in c_0(\R^{2n})$, that is
$$
(-\D)^n w=tK(x) e^{2n\alpha u_0 + 2nw+c_w} - t\alpha(-\Delta)^n u_0 \;\; \mbox{in } \R^{2n}.
$$
Applying again Theorem \ref{linear}, there holds
\begin{align}\label{w}
w(x)=-\frac{t}{\gamma_n} \int_{\R^{2n}} K(y)e^{2n\alpha u_0 + 2nw+c_w} \ln |x-y| dy - t\alpha u_0(x).
\end{align}
Denote
\begin{align}
\label{ww}
\ww := w + \frac{c_w}{2n} + \frac{\ln t}{2n} +t\alpha u_0 \quad \mbox{and} \quad \overline w := \ww - t\alpha u_0 = w + \frac{c_w}{2n} + \frac{\ln t}{2n}.
\end{align}
Clearly, $(-\D)^n\ww= \wK e^{2n\ww}$ in $\R^{2n}$ with $\wK :=K e^{(1-t)2n\alpha u_0}$. As $K \not\equiv 0$, without loss of generality, up to an orthogonal transform of $\R^{2n}$, we can assume that $K < 0$ in $B_3$, so $|\wK|\geq C_0>0$ in $B_3$ for any $t\in (0, 1]$. We deduce then
\begin{align*}
 \int_{B_3} e^{2n\ww} dx\leq \frac{1}{C_0}\int_{B_3} |\wK| e^{2n\ww}dx & = \frac{1}{C_0}\int_{B_3} |K| e^{2nw+c_w+2n\alpha u_0} dx\\
& \leq C_1\int_{\R^{2n}} |K| e^{2nw+c_w+2n\alpha u_0} dx = C_1\alpha \gamma_n.
\end{align*}

\medskip
On the other hand, thanks to \eqref{w}, we have $$\D\ww(x) =-\frac{(2n-2)t}{\gamma_{n}}\int_{\R^{2n}}\frac{\wK e^{2n\ww}}{|x-y|^2}dy.$$
With Fubini's theorem, as $\wK \leq 0$ and $t \in (0, 1]$,  there holds
\begin{align*}
\int_{B_3}|\D\ww| dx \leq \frac{2n-2}{\gamma_n}\int_{\R^{2n}} \wK e^{2n\ww}\int_{B_3(y)}\frac{dx}{|x-y|^2} dy \leq C\alpha.
\end{align*}
Using Theorem 4.2 in \cite{HM}, there exists $C>0$ (independent of $\ww$) such that $$\sup_{B_2}\ww \leq C.$$
From the uniform upper bound of $\ww$, $\overline w$ is uniformly upper bounded in $B_2$.

\medskip
Moreover, we have
\begin{align}
\label{wlocal}
\|\ww\|_{C^{2n-1}(\overline B_1)} + \|\overline w\|_{C^{2n-1}(\overline B_1)}\leq C.
\end{align}
Indeed, let $x\in B_1$, for any $1\leq k \leq 2n-1$, we have
\begin{align*}
|\nabla^k\ww(x)|\leq C_k\int_{\R^{2n}}\frac{|\wK|e^{2n\ww}}{|x-y|^k}dy
& = C_k\int_{B_1(x)}\frac{|\wK|e^{2n\ww}}{|x-y|^k}dy+C_k\int_{\R^{2n}\setminus {B_1(x)}}\frac{|\wK| e^{2n\ww}}{|x-y|^k}dy\\
& \leq C\int_{B_1(x)}\frac{dy}{|x-y|^k} + C\int_{\R^{2n}}|\wK| e^{2n\ww} dy\\
& \leq C + C\alpha\gamma_n.
\end{align*}

\medskip
{\sl Step 4: Uniform upper bound in exterior domain.} For that, we will use the Kelvin transform, let $\wv(x) : = \overline w(x|x|^{-2})$ for $0< |x| \leq 2$, then
\begin{align*}
(-\D)^n \wv = \overline K(x) e^{2n\wv} =: \widetilde f  \;\; \mbox{in } B_2\setminus\{0\},\quad\mbox{where } \overline K(x) = |x|^{-4n-2n\alpha}K\left(\frac{x}{|x|^2}\right).
\end{align*}
Here we used $\D^n u_0(x) \equiv 0$ for $|x| \geq \frac{1}{2}$. As $\Lambda_u > -\alpha_1(K)\gamma_n$, there exists $\e_0>0$ such that $\alpha+\frac{\e_0}{2n}<\alpha_1$, so $K(x)(1+|x|)^{2n\alpha+\e_0}\in L^1(\R^{2n})$, hence $\overline K(x)|x|^{-\e_0}\in L^1(B_2)$. Fix $\e_1 > 0$ such that $(4n+2n\alpha+\lambda)\epsilon_1 = \e_0$ where $\lambda$ is the constant in \eqref{growp}, there holds, similarly as in the proof of Theorem \ref{expan1},
$$
|\overline K(x)|^{1+\e_1}\leq C\overline K(x)|x|^{-\e_0} \left(|x|^{-4n-2n\alpha-\lambda+\frac{\e_0}{\e_1}}\right)^{\e_1} = C\overline K(x)|x|^{-\e_0} \quad \mbox{in } B_2.
$$
Hence $\widetilde f \in L^{1+\e_1}(B_2)$. Let $v_1$ be the solution to
\begin{align*}
(-\D)^n v_1 = \widetilde f \;\; \mbox{in } B_2, \quad v_1 = \D v_1 = \ldots = \D^{n-1} v_1 = 0 \;\; \mbox{on }\p B_2.
\end{align*}
Clearly $v_1 \in W^{2n, 1+ \e_1}(B_2)$. Hence $v_2 = \wv - v_1 + C_t$ verifies
\begin{align*}
\D^n v_2 = 0 \;\; \mbox{in } B_2\setminus\{0\}.
\end{align*}
Here $C_t$ is a large enough constant depending on $t$ such that $v_2 \geq 0$. Pay attention that we do not claim any uniform control on $C_t$ actually, but it does exist for each $t$ since $\overline w$ is bounded at infinity.

\smallskip
Applying Theorem 3.1 in \cite{GMT}, we claim that $v_2$ can be extended as a smooth polyharmonic function in $B_2$, therefore $\wv \in W^{2n, 1+ \e_1}_{loc}(B_2)$ and $(-\D)^n \wv = \widetilde f \leq 0$ in $B_2$.
Denote $\wv_k = (-\D)^k\wv$ for $1 \leq k \leq n-1$. By \eqref{wlocal}, there holds
\begin{align*}
\wv_{n-1} (x) \leq \sup_{\p B_1}\wv_{n-1} \leq C_1\;\; \mbox{in } B_1.
\end{align*}
Consequently,
\begin{align*}
(-\D)\Big(\wv_{n-2} + \frac{C_1}{4n}r^2\Big) \leq 0 \;\; \mbox{in } B_1.
\end{align*}
Hence
\begin{align*}
\wv_{n-2} (x) \leq \sup_{\p B_1}\wv_{n-2} + \frac{C_1}{4n} \leq C_2\;\; \mbox{in } B_1.
\end{align*}
By induction, we can claim finally that $\wv \leq C$ in $B_1$, so we have $\overline w(x) \leq C$ for any $|x| \geq 1$ hence $\overline w\leq C$ in $\R^{2n}$.

\medskip
{\sl Step 5: Conclusion.} Combining {\it Steps 3--4}, if $w\in c_{0}(\R^{2n})$, $w=t\T(w)$  with $t\in(0,1]$, then $2n\overline w = 2nw + c_w + \ln t \leq C$ in $\R^{2n}$.
Recall that $w=t\T(w)$ is equivalent to
\begin{align*}
(-\Delta)^n w = \Psi(w) \;\; \mbox{with } \Psi(w) = K(x)e^{2n\alpha u_0 + 2n\overline w} - t\alpha(-\D)^nu_0.
\end{align*}
Using now the uniform integrability and the uniform polynomial control of $\Psi(w)$, the proof of Theorem \ref{expan0} yields $\|w\|_\infty \leq C$. By Leray-Schauder's fixed point Theorem, there exists a fixed point $w$ for $\T$, so we get the desired solution as $u= w +\alpha u_0 + \frac{c_w}{2n}$. The proof of Theorem \ref{negative} is completed. \qed

\medskip
\section{Normal metrics for radial nonpositive $Q$-curvature}
\reset
In this section, we consider radial functions $K\leq 0$. We will prove that in the radial setting, Theorem \ref{negative} remains true without any growth control on the $Q$-curvature $K$.

\subsection{Radial solution to polyharmonic equation}
Our proof is based on the understanding of linear equation in the radial setting.
\begin{prop}
\label{rlinear0}
Assume that $n \geq 2$ and $\e > 0$, $f \in \dot{L}_{\rm rad}^1(\R^{2n})$ and $f(1+|x|)^\e \in L^1(\R^{2n})$. Then there exists a unique solution $u \in c_{0, {\rm rad}}(\R^{2n})$ to $\D^n u = f$ in $\R^{2n}$. Furthermore, there is $C > 0$ depending only on $n$ and $\e$ such that
\begin{align}
\label{radialv}
\|u\|_\infty \leq C\int_{\R^{2n}} (1+|x|)^\e |f(x)|dx.
\end{align}
\end{prop}

\noindent
{\sl Proof}. As $f \in \dot{L^1}(\R^{2n})$, let
\begin{align}
\label{v}
\begin{split}
v_1(r) & = -\frac{1}{2n-2}\int_r^\infty sf(s)ds -\frac{1}{2n-2}r^{2-2n}\int_0^rs^{2n-1}f(s)ds\\
 & = -\frac{1}{2n-2}\int_r^\infty sf(s)ds +\frac{1}{2n-2}r^{2-2n}\int_r^\infty s^{2n-1}f(s)ds.
\end{split}
\end{align}
Easily $\D v_1 = f$ in $\R^{2n}$. Using the first line of \eqref{v}, we have
\begin{align*}
\| v_1\|_\infty \leq \frac{1}{2n-2} \int^\infty_0 s|f(s)| ds.
\end{align*}
Moreover, by the second line in \eqref{v}, for any $r \geq 1$,
\begin{align*}
| v_1(r)| & \leq \frac{1}{n-1} r^{2-2n} \int^\infty_r s^{2n-1}|f(s)| ds\\
& \leq \frac{1}{(n-1)r^{2n-2+\e}} \int^\infty_r s^{2n-1+ \e}|f(s)| ds\\
& \leq CM_{f,\e}r^{-2n+2-\e}
\end{align*}
where
\begin{align*}
M_{f,\e} := \int_{\R^{2n}} (1+|x|)^\e |f(x)|dx.
\end{align*}
Therefore, there holds
\begin{align}
\label{estv1}
| v_1(r)| \leq C_1M_{f,\e}\Big[\chi_{[0, 1)}(r) + \frac{1}{r^{2n-2+\e}}\chi_{[1, \infty)}(r)\Big], \quad \forall\; r \geq 0.
\end{align}

Now we will apply the following elementary result, which proof is postponed.
\begin{lem}
\label{h}
Let
\begin{align*}
h_\beta(r) = \chi_{[0, 1)}(r) + \chi_{[1, \infty)}(r)\frac{1}{r^\beta}.
\end{align*}
Then for any $2 < \beta < 2n$, there exists a unique solution $w \in c_{0, {\rm rad}}(\R^{2n})$ to $\D w = h_\beta$ in $\R^{2n}$ and we have $|w(r)| \leq C_{n, \beta} h_{\beta-2}(r)$ for all $r \geq 0$.
\end{lem}
Remark that \eqref{estv1} means just $|v_1(r)| \leq C_1M_{f,\e}h_{2n-2+\e}(r)$ in $\R_+$. By induction with Lemma \ref{h}, we get $v_i \in c_{0, {\rm rad}}(\R^{2n})$ such that $\D v_{i+1} = v_i$ in $\R^{2n}$, for $1 \leq i \leq n-1$ and
\begin{align*}
|v_i(r)| \leq C_iM_{f,\e}h_{2n-2i+\e}(r) \;\; \mbox{for all } r \geq 0, \; 1 \leq i \leq n.
\end{align*}
Obviously $u = v_n$ is the desired solution to $\D^n u = f$. \qed

\begin{rem}
\label{eradial}
The above proof yields a more precise estimate at infinity: $|u(r)| \leq Cr^{-\e} M_{f,\e}$ for all $r \geq 1$.
\end{rem}

\medskip
\noindent
{\sl Proof of Lemma \ref{h}}. Let
\begin{align*}
w(r) = -\frac{1}{2n-2}\int_r^\infty h_\beta(s)sds -\frac{1}{2n-2}r^{2-2n}\int_0^rs^{2n-1}h_\beta(s)ds =: -\frac{1}{2n-2}(J_1 + J_2).
\end{align*}
We have, for $r \geq 1$,
$$J_1(r) = \int_r^\infty s^{1- \beta} ds = \frac{1}{\beta - 2}r^{2-\beta}$$
and
\begin{align*}
J_2(r) = r^{2-2n}\int_0^rs^{2n-1}h_\beta(s)ds & = \frac{r^{2-2n}}{2n} + r^{2-2n}\int_1^r s^{2n-1 -\beta} ds\\
& \leq \frac{r^{2-2n}}{2n} +  \frac{1}{2n - \beta}r^{2-\beta}\\
& \leq C r^{2-\beta}.
\end{align*}
For $r \in [0, 1]$, there holds $|w(r)| \leq \|w\|_\infty \leq C\|h_\beta(s)s\|_{L^1(\R_+)}.$  So we are done.\qed

\medskip
By the second line of \eqref{v}, we see that
\begin{align*}
v_1(r) = \frac{1}{2n-2}\int_r^\infty sf(s)\Big[\left(\frac{s}{r}\right)^{2n-2} - 1\Big]ds
\end{align*}
which yields another simpler fact.
\begin{lem}
\label{sign}
Assume that $f \in \dot{L}^1(\R^{2n})$ is radial and $f$ has constant sign for $r \geq r_0 > 0$. Then the unique ground state solution to $\D v = f$ in $\R^{2n}$ has the same sign of $f$ in $\R^{2n}\backslash B_{r_0}$.
\end{lem}

By working more carefully, we get the following result for $f$ having less integrability assumption, which has its own interests.
\begin{thm}
\label{rlinear}
Let $n \geq 2$. Assume that $f \in \dot{L}_{\rm rad}^1(\R^{2n})$, $f(\ln(2+|x|) \in L^1(\R^{2n})$, and $sf(s) \in L^1(\R_+)$. Then there exists a unique $u \in c_{0, {\rm rad}}(\R^{2n})$ solution to $\D^n u = f$ in $\R^{2n}$, and there is $C > 0$ depending only on $n$ such that
\begin{align}
\label{radialv2}
\|u\|_\infty \leq C\int_0^\infty \Big[s + s^{2n-1}\ln(2+|s|)\Big] |f(s)|ds.
\end{align}
\end{thm}

\noindent
Proof. For any radial function $g \in L^1_{loc}(\R^{2n})$, we know that if $sg(s) \in L^1(\R_+)$, then $\Delta \Psi(g) = g$ in $\R^{2n}$ where
\begin{align}
\label{Psi}
\begin{split}
\Psi(g)(r) & := -\frac{1}{2n-2}\int_r^\infty sg(s)ds -\frac{1}{2n-2}r^{2-2n}\int_0^rs^{2n-1}g(s)ds, \quad \forall\; r > 0.
\end{split}
\end{align}

Let $f$ be radial such that $(s + s^{2n-1})f(s) \in L^1(\R_+)$. We define for any $1 \leq \ell \leq 2n-1$, $k \in \R$ and $r > 0$,
\begin{align*}
H_{k, \ell}(r) = r^k\int_r^\infty s^\ell f(s) ds, \quad \mbox{and}\quad G_{k, \ell}(r) = r^k\int_0^r s^\ell f(s) ds.
\end{align*}
The rest of proof is based on the following lemma which can be proved by direct calculations.
\begin{lem}
\label{H}
Let $(s + s^{2n-1})f(s) \in L^1(\R_+)$.
\begin{itemize}
\item[(i)] If $k > -2n$, $1\leq \ell \leq 2n- 1$, $1\leq k + \ell + 2 \leq 2n-1$, there holds $H_{k, \ell} \in L^1_{loc}(\R^{2n})$, $sH_{k, \ell}(s) \in L^1(\R_+)$ and
\begin{align*}
\Psi(H_{k, \ell}) \in {\rm Vect}\Big\{ H_{0, k+\ell+2};\; H_{k+2, \ell};\; G_{2-2n, k+\ell+2n}\Big\}.
\end{align*}
\item[(ii)] If $-2n < k < -2$, $1 \leq k + \ell + 2 \leq 2n-1$, there holds $G_{k, \ell} \in L^1_{loc}(\R^{2n})$, $sG_{k, \ell}(s) \in L^1(\R_+)$ and
\begin{align*}
\Psi(G_{k, \ell}) \in {\rm Vect}\Big\{ H_{0, k+\ell+2};\; G_{k+2, \ell};\; G_{2-2n, k+\ell+2n}\Big\}.
\end{align*}
\end{itemize}
\end{lem}

Consider now a continuous radial function $f \in L^1(\R^{2n}, \ln(2+|x|) dx )\cap \dot{L}^1(\R^{2n})$. Define $v_1 = \Psi(f)$, then
\begin{align*}
v_1 =  -\frac{1}{2n-2}H_{0, 1} -\frac{1}{2n-2}G_{2-2n, 2n-1} = -\frac{1}{2n-2}H_{0, 1} + \frac{1}{2n-2}H_{2-2n, 2n-1}.
\end{align*}
Using Lemmas \ref{H} (i), we have $v_2 = \Psi(v_1)$ is well defined if $n \geq 3$ and
\begin{align*}
v_2 \in {\rm Vect}\Big\{ H_{0, 3};\; H_{2, 1};\; H_{4-2n, 2n-1}; \; G_{2-2n, 1+2n}\Big\}.
\end{align*}
More generally, let $v_k = \Psi(v_{k-1})$ for $2 \leq k \leq n-1$, we can check by induction, using Lemma \ref{H} that they are well defined, and for any $1 \leq k \leq n-1$, $v_k$ belongs to
\begin{align*}
{\rm Vect}\Big\{ H_{2i, 2k-2i-1};\; i \leq k-1\Big\} \oplus{\rm Vect}\Big\{ H_{2k-2n, 2n-1}\Big\} \oplus {\rm Vect}\Big\{G_{2i-2n, 2n+2k-2i-1}; \; 1 \leq i \leq k-1\Big\}.
\end{align*}
In particular,
\begin{align*}
v_{n-1} = a_n H_{-2, 2n-1} + \sum_{i = 0}^{n-1} b_{i, n} H_{2i, 2n-2i-3} + \sum_{i = 1}^{n-2}c_{i, n} G_{2i-2n, 4n-2i-3} =: a_n H_{-2, 2n-1} +  J.
\end{align*}
Here $a_n$, $b_{i, n}$ and $c_{i, n}$ are universal constants depending on $i$ and $n$. Applying again Lemma \ref{H},
$$\int_0^{+\infty} s|J(s)| ds \leq C\int_0^\infty \left(s + s^{2n-1}\right) |f(s)| ds.$$

It remains to study $H_{-2, 2n-1}$. Recall that $f \in \dot{L}^1(\R^{2n})$, so $H_{-2, 2n-1} = -G_{-2, 2n-1}$. Therefore
\begin{align*}
\int_0^\infty s|H_{-2, 2n-1}(s)|ds & = \int_0^1 s|G_{-2, 2n-1}(s)|ds + \int_1^\infty s|H_{-2, 2n-1}(s)|ds\\
& \leq \int_0^1\frac{1}{s}\int_0^s \sigma^{2n-1}|f(\sigma)|d\sigma + \int_1^\infty\frac{1}{s}\int_s^\infty \sigma^{2n-1}|f(\sigma)|d\sigma\\
& = \int_0^\infty\sigma^{2n-1}|f(\sigma)||\ln\sigma| d\sigma\\
& \leq C\int_0^1\sigma|f(\sigma)| d\sigma + C\int_1^\infty \sigma^{2n-1}|f(\sigma)||\ln\sigma| d\sigma.
\end{align*}
To conclude, $sv_{n-1}(s) \in L^1(\R_+)$ and $v_{n-1} \in L^1_{loc}(\R^{2n})$, there exists then a unique $u \in c_{0, {\rm rad}}(\R^{2n})$ such that $\Delta u = v_{n-1}$ in $\R^{2n}$. Moreover, $\|u\|_\infty \leq C\|sv_{n-1}(s)\|_{L^1(\R_+)}$, so we are done. \qed

\subsection{Proof of Theorem \ref{K<0radial}} Fix $u_0 \in C^\infty_{\rm rad}(\R^{2n})$ such that $u_0(r) = \ln r$ for $r \geq 1$. Clearly, $\Delta^n u_0\equiv0$ in $\R^{2n}\backslash B_1$ and
$$\int_{\R^{2n}} (-\D)^n u_0 dx = -\gamma_n.
$$
Define $w = u + \alpha u_0.$ Then $u$ is a solution of \eqref{0} if and only if $w$ satisfies $$\Delta^n w = (-1)^nK e^{2n\alpha u_0 + 2nw} - \alpha \Delta^n u_0.$$ To construct such $w$, we will use again Leray-Schauder's fixed point argument as before.

\medskip
Let $\e > 0$, consider the space
\begin{align*}
E_\e = L^1(\R^{2n}, (1+|x|)^\e dx)\cap \dot{L}^1(\R^{2n}),
\end{align*}
endowed the norm $\|(1+|x|)^\e f\|_{L^1(\R^{2n})}$. Define the operator $$\T: E_\e \longrightarrow \left(c_{0,{\rm rad}}(\R^{2n}), \|\cdot\|_{L^\infty(\R^{2n})}\right)$$ as follows. For any $f \in E_\e$, $\T(f) = w$ is the unique solution in $c_{0,{\rm rad}}(\R^{2n})$ such that $\D^n w = f$ in $\R^{2n}$ given by Proposition \ref{rlinear0}. With Remark \ref{eradial}, we can prove easily that $\T$ is compact.

\medskip
Let $\alpha \in (0, \alpha_1)$. As above, for any $w \in c_{0,{\rm rad}}(\R^{2n})$, there is a unique constant $c_w \in \R$ such that $(-1)^n Ke^{2n\alpha u_0+2nw+c_w} - \alpha\D^nu_0 \in \dot{L}^1(\R^{2n})$. Fix $0 < \e <\alpha_1-\alpha$, so $(1+|x|)^{\e} Ke^{2n\alpha u_0}\in L^1(\R^{2n})$. Define a mapping $\Phi$ from $c_{0,{\rm rad}}(\R^{2n})$
into itself by
\begin{align*}
\Phi(w) = \T\Big( (-1)^n Ke^{2n\alpha u_0+2nw+c_w} - \alpha\D^n u_0\Big).\end{align*}
Clearly, $\Phi$ is compact. Now we need only to check that the family of eventual fixed points of $t\Phi$ with $t\in (0, 1]$ is uniformly bounded. Assume that $t\Phi(w) = w$, $t \in (0, 1]$.

\medskip
If $K$ is compactly supported. Thanks to Proposition \ref{rlinear0}, there holds
\begin{align*}
\|w\|_{L^\infty(\R^{2n})} & \leq Ct\int_{\R^{2n}} (1+|x|)^\e \left|Ke^{2n\alpha u_0+2nw+c_w} - (-1)^n\alpha\D^nu_0\right| dx\\
& \leq C\int_{\R^{2n}} |K|e^{2n\alpha u_0+2nw+c_w}dx + C =  C\gamma_n + C.
\end{align*}
Hence we get the uniform estimate of $w$.

\medskip
Suppose now ${\rm supp}(K)$ is unbounded, we fix $R_0 > 1$ and $\eta > 0$ such that $K < 0$ in $[R_0, R_0+\eta]$. Define $A:=\overline{B}_{R_0+\eta}\setminus B_{R_0}$ and $C_0 := \max_A K < 0$. Let
$$\ww = w+ \frac{c_w}{2n} + \frac{\ln t}{2n} + t\alpha u_0.$$
Then
$(-\D)^n\ww = Ke^{2n(1-t)\alpha u_0}e^{2n\ww} \leq 0$ in $\R^{2n}$, we denote $\wv_k = (-\Delta)^{n-k}\ww$.

\medskip
Notice that $\D \wv_1 \geq 0$ in $\R^{2n}$, which means $\wv_1$ is increasing in $r$. As $\ww = O(\ln r)$ at $\infty$, there must holds $\lim_{r\to \infty}\wv_1(r) =0$, hence $\wv_1 \leq 0$ in $\R^{2n}$. By induction, we deduce that for any $1 \leq k \leq n-1$, $\wv_k$ is negative, increasing in $r$ and tends to $0$ at infinity. So $\D \ww = - \wv_{n-1} \geq 0$ in $\R^{2n}$, and then $\ww$ is nondecreasing in $r$. Therefore,
\begin{align*}
|C_0| e^{2n(1-t)\alpha \ln R_0}e^{2n\ww(R_0)}|A| & \leq \int_A |K|e^{2n(1-t)\alpha u_0}e^{2n\ww} dx\\
& \leq \|Ke^{2n(1-t)\alpha u_0}e^{2n\ww}\|_{L^1(\R^{2n})} = \gamma_n\alpha t.
\end{align*}
We deduce then
$$
e^{2n\ww(R_0)}\leq \frac{8\gamma_n\alpha t}{|C_{R_0}| |A|}e^{2n(t-1)\alpha \ln R_0}\leq C,
$$
where $C$ is independent of $t\in(0,1]$. In other words, $\max_{B_{R_0}}\ww(r) = \ww(R_0) \leq C'.$

\smallskip
Consequently, there is $C_1 > 0$ independent on $t\in(0,1]$ such that
\begin{equation}\label{2.9}
\sup_{B_1} \hw \leq C_1, \quad \mbox{where }\; \hw :=  w + \frac{c_w}{2n} + \frac{\ln t}{2n}.
\end{equation}
As before, we can claim that $\|\nabla^k \hw\|_{L^\infty(\overline B_1)} \leq C$ for all $0 \leq k \leq 2n - 1$.

\medskip
On the other hand, denote $\hv_k = (-\Delta)^{n-k}\hw$ for $1 \leq k \leq n-1$. Remark first that
\begin{align}
\label{wvk}
\forall\; 1 \leq k \leq n-1,\quad \hv_k(r) = \wv_k(r) - t\alpha(-\D)^{n-k}u_0 \rightarrow 0\;\; \mbox{as } r \to \infty.
\end{align}
Recall also ${\rm supp}(\D^n u_0) \subset B_1$, then for $r \geq 1$,
$-\D \hv_1 =(-\D)^n \ww = Ke^{2n\alpha u_0+2n\hw} \leq 0$. By Lemma \ref{sign}, we get $\hv_1 \geq 0$ in $\R^{2n}\backslash B_1$.

\medskip
As $-\D \hv_2 = \hv_1$, we get $\D \hv_2 \leq 0$ if $r \geq 1$ so that $r^{2n-1}{\hv_2}'(r) \leq \hv_2'(1)$ for $r \geq 1$, i.e.~$\hv_2'(r) \leq Cr^{1-2n}$ thanks to the uniform estimate in $\overline B_1$. Integrating from $r$ to $\infty$ and using \eqref{wvk}, there holds $\hv_2'(r) \geq -Cr^{2-2n}$ for any $r \geq 1$.

\medskip
By induction, we can claim that for any $1 \leq k \leq n-1$, there is a universal positive constant $C_k$ such that
\begin{align*}
\hv_k(r) \geq -C_k r^{2k-2-2n} \quad \forall\; r \geq 1.
\end{align*}
In particular, let $k = n-1$, we have $\D\hw \leq C_{n-1}r^{-4}$ for $r \geq 1$. Therefore, for $r \geq 1$
\begin{align*}
r^{2n-1}\hw'(r) = \int_0^r s^{2n-1} \D\hw(s) ds \leq C + C\int_1^r s^{2n-5} ds \leq C + Cr^{2n-3}.
\end{align*}
Hence $\hw'(r) \leq Cr^{-2}$ for $r \geq 1$, which yields $\hw(r) \leq C$ for $r \geq 1$.

\medskip
To conclude, $\hw$ is uniformly upper bounded in $\R^{2n}$, so $e^{2n\widehat w}$ is uniformly bounded in $\R^{2n}$. It's easy to see that $\|w\|_\infty \leq C$, as
$$w = t\Phi(w) = \T\Big(Ke^{2n\alpha u_0 + 2n\hw} - t\alpha(-\D)^n u_0\Big).$$
We can apply for example \eqref{radialv}. By Leray-Schauder's argument, there exists a fixed point $w$ for $\Phi$, hence the desired solution. \qed

\bigskip
\noindent
{\bf Acknowledgements.} X. Huang and D. Ye are partially supported by NSFC (No.~12271164). F. Zhou is partially supported by NSFC (No.~12071189). All of them are partially supported by Science and Technology Commission of Shanghai Municipality (No.~22DZ2229014).


\medskip


\begin{thebibliography}{99}

\bibitem{BO}
T.P. Branson and B. Orsted, Explicit functional determinants in four dimensions, Proc. AMS. 113 (1991) 669-682.


\bibitem{BHL}
F. Brito, J. Hounie and M.L. Leite, Liouville's formula in arbitrary planar domains, Nonlinear Anal., 60(7) (2005), 1287-1302.

\bibitem{CC}
S.-Y.A. Chang and W. Chen, A note on a class of highter order conformally covariant equations, Discrete Contin. Dyn. Syst. 7(2) (2001), 275-281.


\bibitem{CQY}
S.-Y.A. Chang, J. Qing and P.C. Yang, On the Chern-Gauss-Bonnet integral for conformal metrics on $\R^4$, Duke Math. J. 103(3) (2000), 523-544.

\bibitem{CL}
W. Chen and C. Li, Classification of solutions of some nonlinear elliptic equations, Duke Math. J. 63(3) (1991), 615-622.

\bibitem{CL1}
K.S. Cheng and C.S. Lin, On the asymptotic behavior of solutions of the conformal Gaussian curvature equations in $\mathbb{R}^2$, Math. Ann. 308(1) (1997), 119-139.


\bibitem{CL3}
 K.S. Cheng and C.S. Lin, Conformal metrics in $\mathbb{R}^2$ with prescribed Gaussian curvature with positive total curvature, Nonlinear Anal. 38(6) (1999), Ser. A: Theory Methods, 775-783.

\bibitem{CL4}
K.S. Cheng and C.S. Lin, Multiple solutions of conformal metrics with negative total curvature, Adv. Differential Equations 5(10-12) (2000), 1253-1288.

\bibitem{CN1}
K.S. Cheng and W-.M. Ni, On the structure of the conformal Gaussian curvature equation on $\mathbb{R}^2$, Duke Math. J. 62(3) (1991), 721-737.

\bibitem{CN2}
K.S. Cheng and W-.M. Ni, On the structure of the conformal Gaussian curvature equation on $\mathbb{R}^2$. II,  Math. Ann. 290(4) (1991), 671-680.

\bibitem{CYZ} H.Y. Chen, D. Ye and F. Zhou, On Gaussian curvature equation in $\R^2$  with prescribed nonpositive curvature, Discrete Contin. Dyn. Syst. 40(6) (2020), 3201-3214.

\bibitem{CW}
K.S. Chou and T. Wan, Asymptotic radial symmetry for solutions of $\Delta u + e^u =0$ in a punctured disc, Pacific J. Math., 1653(2) (1994), 269-276.

\bibitem{GMT} M. Ghergu, A. Moradifam and S.D. Taliaferro, Isolated singularities of polyharmonic inequalities, J. Funct. Anal. 261(3) (2011), 660-680.

\bibitem{GT}
D. Gilbarg and N. Trudinger, Elliptic partial differential equations of second order, Reprint of the 1998 edition. Classics in Mathematics. Springer-Verlag, Berlin, 2001.


\bibitem{GJMS}
C. Graham, R. Jenne, L. Mason and G. Sparling, Conformally invariant powers of the Laplacian. I: existence, J. Lond. Math. Soc., II. Ser. 46(3) (1992), 557-565.

\bibitem{GHYZ}
Z.M. Guo, X. Huang, D. Ye and F. Zhou, Qualitative properties of H\'enon type equations with exponential nonlinearity, Nonlinearity 35(1) (2022), 492-512.

\bibitem{Gu}
M.J. Gursky, The principal eigenvlue of a conformally invariant differential operator, with an application to semlinear elliptic PDE, Comm. Math. Phys. 207 (1999), 131-143.

\bibitem{HY}
X. Huang and D. Ye, Conformal metrics in $\R^{2m}$ with constant $Q$-curvature and arbitrary volume, Calc. Var. Partial Differential Equations 54(4) (2015), 3373-3384.

\bibitem{Hu}
A. Huber, On subharmonic functions and differential geometry in the large, Comment. Math. Helv. 32 (1957) 13-72.

\bibitem{H1}
A. Hyder, Existence of entire solutions to a fractional Liouville equation in $\R^n$, Rend. Lincei Mat. Appl. 27 (2016), 1-14.

\bibitem{H}
A. Hyder, Conformally Euclidean metrics on $\R^n$ with arbitrary total $Q$-curvature, Anal. PDE 10(3) (2017), 635-652.

\bibitem{HMM}
A. Hyder, G. Mancini and L. Martinazzi, Local and nonlocal singular Liouville equations in Euclidean spaces,
Int. Math. Res. Not. 15 (2021), 11393-11425.

\bibitem{HM}
A. Hyder and L. Martinazzi, Conformal metrics on $\mathbb{R}^{2m}$ with constant $Q$-curvature, prescribed volume and asymptotic behavior, Discrete Contin. Dyn. Syst. 35(1) (2015), 283-299.

\bibitem{HM1}
A. Hyder and L. Martinazzi, Normal conformal metrics on $\R^4$ with $Q$-curvature having power-like growth, J. Differential Equations 301 (2021), 37-72.

\bibitem{JSTW}
J. Jin, L. Shu, D. Tai and D. Wu, Normal conformal metrics on $\R^4$ with singular $Q$-curvature having power-like growth, J. Differential Equations 327 (2022), 64-108.

\bibitem{Li1}
M. Li, The total $Q$-curvature, volume entropy and polynomial growth polyharmonic functions, Adv. Math. 450 (2024), Paper no. 109768, 43 pp.

\bibitem{Lin}
 C.S. Lin, A classification of soluions of conformally invariant fourth order equations in $\R^n$, Comment. Math. Helv. 73 (1998), 206-231.

\bibitem{L1}
J. Liouville, Sur le th\'eor\`eme de M. Gauss, concernant le produit des deux rayons de courbure
principaux en chaque point d'une surface, Note IV to G. Monge //  Applications de Analyse
\`a la G\'eom\'etrie / Paris: Bachelier. (1850), 583-600.

\bibitem{L2}
J. Liouville, Sur  l'\'equation  aux  diff\'erences  partielles  $\frac{d^2\log \lambda}{du dv} \pm \frac{\lambda}{2a^2}=0$, Journal de
Math\'ematiques Pures et Appliqu\'ees, 18 (1853), 71-72.

\bibitem{M2}
L. Martinazzi, Conformal metrics on $\R^{2m}$ with constant $Q-$curvature, Rend. Lincei. Mat.
Appl. 19 (2008), 279-292.

\bibitem{M}
L. Martinazzi, Conformal metrics on $\mathbb{R}^{2m}$ with constant $Q$-curvature and large volume, Ann. I.H.P. Analyse non lin\'eaire 30 (2013), 969-982.

\bibitem{Mc1}
R.C. McOwen, On the equation $\Delta u +K e^{2u}=f$ and prescribed negative curvature in $\mathbb{R}^2$, J. Math. Anal. Appl. 103(2) (1984), 365-370.

\bibitem{Mc2}
R.C. McOwen, Conformal metrics in ${\mathbb R}^2$ with prescribed Gaussian curvature and positive total curvature, Indiana Univ. Math. J. 34(1) (1985), 97-104.


\bibitem{Ni}
W-.M. Ni, On the elliptic equation $\Delta u + K(x)e^{2u}=0$ and conformal metrics with prescribed Gaussian curvatures, Invent. Math. 66 (1982), 343-352.

\bibitem{O}
R. Osserman, On the inequality $\Delta u\geq f(u)$, Pacific J. Math. 7 (1957), 1641-1647.

\bibitem{P}
S. Paneitz, A quartic conformally covariant differential operator for arbitrary pseudo-Riemannian manifolds, (Preprint, 1983), SIGMA Symmetry Integrability Geom. Methods Appl. 4 (2008), 036, 3 pages.

\bibitem{WW}
S. Wang and Y. Wang, Integrability of scalar curvature and normal metric on conformally flat manifolds, J. Differential Equations 265(4) (2018), 1353-1370.

\bibitem{WYM}
Y.M. Wang, Nonexistence of solutions to singular $Q$-curvature equations, Arch. Math. 117(4) (2021), 455-467.

\bibitem{WX} J. Wei and X. Xu, Classification of solutions of higher order conformally invariant equations, Math. Ann. 313 (1999), 207-228.

\bibitem{WY}
J. Wei and D. Ye, Nonradial solutions for a conformally invariant fourth order equation in $\mathbb{R}^4$, Calc. Var. P.D.E. 32(3) (2008), 373-386.
\end{thebibliography}
\end{document}